[1994/12/01]
\documentclass{amsart}
\chardef\bslash=`\\ 

\usepackage{amsmath,amsfonts,amssymb}
\usepackage{mathabx}
\usepackage{graphicx}
\usepackage{url}
\usepackage{hyperref}
\usepackage{comment}
\usepackage{makecell}
\usepackage{boldline}
\setcellgapes{3pt}
\usepackage{hhline,colortbl}
\usepackage{xcolor}
\usepackage{textcomp}

\definecolor{blueviolet}{rgb}{0.5,0,1}

\newtheorem{thm}{Theorem}[section]

\newtheorem{lem}[thm]{Lemma}

\newtheorem{con}[thm]{Conjecture}
\newtheorem{ex}[thm]{Example}

\theoremstyle{definition}

\theoremstyle{remark}

\newcommand{\Z}{\mathbb{Z}}

\UseRawInputEncoding

\newcommand{\eval}[2][\right]{\relax
  \ifx#1\right\relax \left.\fi#2#1\rvert}

\begin{document}
\title[Maximal dilatation]{Maximal dilatation on nonorientable surfaces\thanks{We thank Erwan Lanneau, Livio Liechti, Julien Marché, Alan Reid, and Darren Long.}}
\author[J. Ham]{Ji-Young Ham}
\address{Department of Mathematics, Kangwon National University, 1 Kangwondaehakgil, Chuncheon-si, Gangwon state, 24341\\
School of Liberal Arts, Seoul National University of Science and Technology, 232 Gongneung-ro, Nowon-gu, Seoul, 01811\\
  Korea} 
\email{jiyoungham@gmail.com}\thanks{The author is supported by Basic Science Research Program through the National Research Foundation of Korea (NRF) funded by the Ministry of Education, Science and Technology (No. NRF-2018005847).}

\author[J. Lee]{Joongul Lee}
\address{Department of Mathematics Education, Hongik University\\
94 Wausan-ro, Mapo-gu, Seoul 04066\\
   Korea}
\email{jglee@hongik.ac.kr\thanks{The author was supported by 2017 Hongik University Research Fund.}} 

\subjclass[2010]{37E30, 37B40, 57M60}

\keywords{Maximal degree, Nonorientable surface, Liechti-Strenner's polynomial, pseudo-Anosov, Minimal dilatation}

\begin{abstract}
 On each nonorientable surface of even genus $g \geq 4$, we show that the Liechti-Strenner's polynomial in~\cite{LiechtiStrenner18} gives a maximal dilatation among pseudo-Anosov diffeomorphisms with an orientable invariant foliation. This is proved by showing that this polynomial is irreducible.
 \end{abstract}
\maketitle

\section{Introduction}
In 1970's, Thurston classified the mapping class group of a surface into periodic, pseudo-Anosov, and reducible~\cite{Thurston88}. The same content from a somewhat different point of view can be found in~\cite{CassonBleiler88}. In this paper, we adopt the Penner's approach~\cite{Penner88,Penner91}. Penner made use of bigon tracks, a slight generalization of train track. Nice example of bigon tracks can be found in~\cite{LiechtiStrenner18-1}. Our paper is based on 
Liechti-Strenner's pseudo-Anosov diffeomorphisms on nonorientable surfaces~\cite{LiechtiStrenner18}.

Let $\Sigma_g$ be a surface of finite type.
A diffeomorphism $h$ of $\Sigma_g$ is called \emph{pseudo-Anosov} if there is a pair of transversely measured foliations 
$\mathcal{F}^u$ and $\mathcal{F}^s$ in $\Sigma_g$ and a real number $\lambda >1$ such that $h(\mathcal{F}^u)=\lambda \mathcal{F}^u$ and  $h(\mathcal{F}^s)=1/\lambda \mathcal{F}^s$~\cite{Thurston88,CassonBleiler88}.
The number $\lambda$ is called the \emph{dilatation} of $h$ and the logarithm of $\lambda$ is called the \emph{topological entropy}. A beginner-friendly example about the dilatation can be found in~\cite{BaikRafiqiWu16}. The set of dilatations of pseudo-Anosov diffeomorphisms of the group of isotopy classes $\Sigma_g$ is discrete~\cite{ArnouxYoccoz81, Ivanov88}. In particular, there exists the minimal dilatation. Note that $\lambda$ is an algebraic number and hence there exists the minimal polynomial of it. Denote by the \emph{degree} of $\lambda$ the degree of the minimal polynomial of it. Denote by a \emph{maximal dilatation} on a set $S$ of diffeomorphisms of 
$\Sigma_g$ a dilatation that gives the maximal degree among $S$.


 Liechti and Strenner~\cite{LiechtiStrenner18} determined the minimal dilatation of pseudo-Anosov diffeomorphisms with an orientable invariant foliation on the closed nonorientable surfaces of genus 4, 5, 6, 7, 8, 10, 12, 14, 16, 18 and 20. 
 Denote by $N_g$ the closed nonorientable surface of genus $g$ and by 
$\delta^{+} (N_g)$ the minimal dilatation among pseudo-Anosov diffeomorphisms of $N_g$ with an orientable invariant foliation. Liechti and Strenner conjectured that the largest root of $x^{2k-1}-x^k-x^{k-1}-1$ is the minimal dilatation of $N_{2k}$ for $k \geq 2$(Conjecture~\ref{con:1}). Denote by \emph{Liechti-Strenner’s polynomials} the odd degree polynomials in the conjecture~\ref{con:1}.
 
\begin{con}~\cite[Conjecture 1.2]{LiechtiStrenner18} \label{con:1}
For all $k \geq 2$, $\delta^{+} (N_{2k})$ is the largest root of $x^{2k-1}-x^k-x^{k-1}-1$.
\end{con}

In~\cite{HamLee20}, we showed that the Liechti-Strenner's example for the closed nonorientable surface $N_{2k}$ in~\cite{LiechtiStrenner18} minimizes the dilatation within the class of pseudo-Anosov diffeomorphisms with an orientable invariant foliation and all but the first coefficient of the characteristic polynomial of the action induced on the first cohomology nonpositive.

The main purpose of this paper is the following two theorems.

\begin{thm} \label{thm:irr}
Denote by $N_{2k}$ the closed nonorientable surface of genus $2k$. 
For all $k \geq 2$, the Liechti-Strenner\textquotesingle s polynomial,
$$ x^{2k-1}-x^k-x^{k-1}-1,$$ is irreducible.
\end{thm}

\begin{thm} \label{thm:main}
For all $k \geq 2$, the Liechti-Strenner\textquotesingle s polynomial,
$$ x^{2k-1}-x^k-x^{k-1}-1,$$ gives a maximal dilatation on $N_{2k}$ among pseudo-Anosov diffeomorphisms of $N_{2k}$ with an orientable invariant foliation.
\end{thm}

 It is known that this degree is possible~\cite{Strenner17}, but there have been no concrete examples of this being realized. The well known conjecture is that this polynomial gives the minimal dilatation~\cite{LiechtiStrenner18}. Here we prove the maximality instead of the minimailty. On each orientable surface of genus $g \geq 2$, a maximal dilatation among pseudo-Anosov diffeomorphisms with orientable invariant foliations is concretely presented in~\cite{LanneauLiechti24,Shin16}, but it doesn't seem to be the minimal dilatation. Earlier work concerning the degree of dilatations can be found in~\cite{ArnouxYoccoz81} and ~\cite{Long85}.
Some references related to minimal dilatations are ~\cite{Abikoff80,AaberDunfield10,Bauer92,BestvinaHandel95,Birman74,Brinkmann00,Brinkmann04,ChoHam08,FarbMargalit12,FarberReinosoWang24,FLP79,Ham06,HamSong07,Hironaka10,HironakaKin06,KinTakasawa13,HironakaTsang24,Leininger04,LanneauThiffeault11,LanneauThiffeault11b,LiechtiStrenner18,LiechtiStrenner21,ErwanLivioCheuk25,Loving19,McMullen15,PapadopoulosPenner87,SnapPy,SongKoLos02,Tsai09,TsangZeng24,Valdivia12,Venzke08,Yazdi20,Zhirov95}.
The above references are neither exhaustive nor systematically assembled, so that their inclusion or noninclusion should not be taken as an indication of quality or relevance of the content. Nonetheless, we hope that you find them helpful.

\section{Liechti-Strenner construction of nonorientable surfaces} \label{section:LS}

We will briefly introduce the Liechti-Strenner's method of constructing the nonorientable surface $\Sigma_{2k-1,2}$ of genus $g=2k$ and the Liechti-Strenner's polynomial.

\subsection{The graph $G_{2k-1,2}$}
\begin{figure}
\begin{center}
\includegraphics[scale=1]{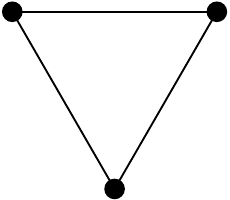}
\caption{The graph $G_{3,2}$.}
\label{fig:graph}
\end{center}
\end{figure}

Let $k \geq 2$ be an odd natural number. Let $G_{2k-1,2}$ be the graph whose vertices are the vertices of a regular $(2k-1)$-gon and every vertex $v$ is connected to the $2$ vertices that are the farthest away from $v$ in the cyclic order of the vertices. Figure~\ref{fig:graph} shows the graph 
$G_{3,2}$.

\subsection{The surface $\Sigma_{2k-1,2}$}
\begin{figure}
\begin{center}
\includegraphics[scale=1]{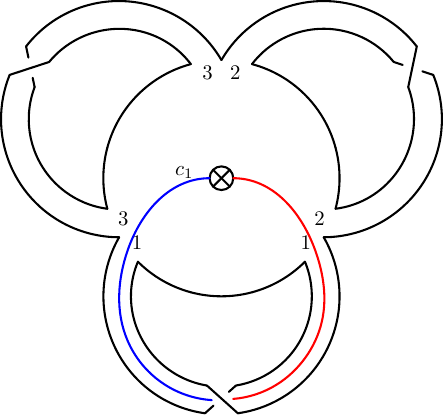}
\caption{The surface $\Sigma_{3,2}$ and the curve $c_1$.}
\label{fig:surface}
\end{center}
\end{figure}

For each $G_{2k-1,2}$, Liechti and Strenner constructed an nonorientable surface $\Sigma_{2k-1,2}$ that contains a collection of curves with
intersection graph $G_{2k-1,2}$~\cite[Subsection 2.2]{LiechtiStrenner18}. 
To construct $\Sigma_{2k-1,2}$, start with a disk with one crosscap. Next, we consider $4k-2$ disjoint intervals on the boundary of the disk and label
the intervals with integers from $1$ to $2k-1$ so that each label is used exactly twice.
In the cyclic order, the labels are $1,s,2,s+1,...,2k,s+2k$ where $s=k+2$ and
all labels are understood modulo $2k-1$.
For each label, the corresponding two intervals are connected by a twisted strip. Figure~\ref{fig:surface} shows the surface 
$\Sigma_{3,2}$.

\begin{lem}~\cite[Proposition 2.3]{LiechtiStrenner18}
The surface $\Sigma_{2k-1,2}$ is homeomorphic to the nonorientable surface of genus $2k$ with
$1$ boundary components.
\end{lem}

\subsection{The curves}
\begin{figure}
\begin{center}
\includegraphics[scale=1]{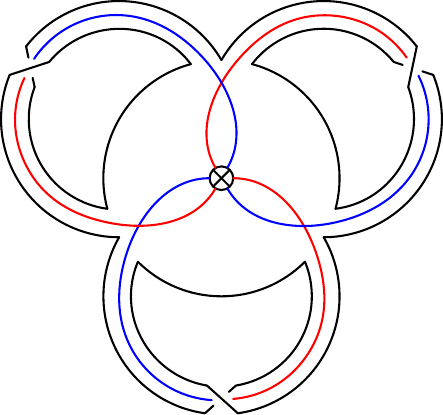}
\caption{A collection of filling inconsistently marked curves.}
\label{fig:curve}
\end{center}
\end{figure}

Liechti and Strenner constructed a two-sided curve $c_i$ for each label $i=1,\ldots,2k-1$ as follows. Each curve consists of two parts.
One part of each curve is the core of the strip corresponding to the label. The other part is an arc inside the disk that passes through the crosscap and connects the corresponding two intervals. The curve $c_1$ is shown in Figure~\ref{fig:surface}.
We choose markings for the $c_i$ which are invariant under the rotational symmetry See Figure~\ref{fig:curve}.

Note that every pair of curves intersects either once or not at all. The curves $c_i$ and $c_j$ are disjoint if and only if the two $i$ labels and the two $j$ labels link in the cyclic order. In other words, if the two $i$ labels separate the two $j$ labels.

\section {The mapping classes}

Denote by $r$ the rotation of $\Sigma_{2k-1,2}$ by one click in the clockwise direction. Denote by $T_{c_1}$ the right-handed Dehn twist about the curve $c_1$.
Define the mapping class
$$\Phi_k=r \circ T_{c_1}.$$
Note that $\Phi_k$ is pseudo-Anosov.

According to Penner~\cite[Theorem 3.1, Theorem 4.1]{Penner88} and Liechti and Strenner~\cite[Proposition 2.6]{LiechtiStrenner18}, we can construct a bigon track on 
$\Sigma_{2k-1,2}$ with filling curves $c_i$, $i=1,\ldots, 2k-1$. Each $c_i$ defines a characteristic measure $\mu_i$ on this bigon track, defined by assigning $1$ to the branches traversed by $c_i$ and zero to the rest. Let $H$ denote the cone generated by the measures $\mu_i$ in the cone of measures on the bigon track. 
Let $\mathcal{C}_k$ be $\mathcal{C}_k=\{c_i\,|i=1,\ldots,2k-1\}$ and $\mathcal{S}_k$ be the semigroup with presentation
$$\mathcal{S}_k (\mathcal{C}_k)=\langle c_i \in \mathcal{C}_k \, : \, c_i \leftrightarrow c_j \text{ if } c_i \cap c_j = \emptyset \rangle.$$

\vspace{2ex}

\begin{thm}~\cite[Theorem 3.4]{Penner88}
The action of $\mathcal{S}_k (\mathcal{C}_k)$ on $H$ admits a faithful representation as a semigroup of invertible (over $\Z$) positive matrices.
\end{thm}

Note that each $c_1$ corresponds to $T_{c_1}$ and the action of $c_1$ on $H$ in the basis $\mu_i$ is given by
$$I_{2k-1}+A$$
where $I_{2k-1}$ is the $(2k-1) \times (2k-1)$-identity matrix, and $A=\begin{bmatrix} a_{ij} \end{bmatrix}$ with $a_{ij}=0$ if $i \neq 1$, and $a_{1j}=\text{card}(c_1 \cap c_j)$.
The rotation $r$ acts by a permutation matrix.
The product of these two matrices is the companion matrix of the Liechti-Strenner's polynomial.

The following example shows how we obtain the Liechti-Strenner's polynomial on $\Sigma_{5,2}$.

\begin{ex}
On $\Sigma_{5,2}$, the action of $c_1$ on $H$ in the basis $\{\mu_i\}$ is given by

$$\begin{bmatrix}
1 & 0 & 1 & 1 & 0  \\
0 & 1 & 0 & 0 & 0  \\
0 & 0 & 1 & 0 & 0  \\
0 & 0 & 0 & 1 & 0  \\
0 & 0 & 0 & 0 & 1  \\

\end{bmatrix}.$$

The action of $r$ on $H$ in the basis $\mu_i$ is given by

$$\begin{bmatrix}
0 & 1 & 0 & 0 & 0 \\
0 & 0 & 1 & 0 & 0  \\
0 & 0 & 0 & 1 & 0 \\
0 & 0 & 0 & 0 & 1 \\
1 & 0 & 0 & 0 & 0  \\
\end{bmatrix}.$$

Hence $\Phi_3$ is represented by the below matrix $M_3$:

$$\begin{bmatrix}
0 & 1 & 0 & 0 & 0  \\
0 & 0 & 1 & 0 & 0  \\
0 & 0 & 0 & 1 & 0  \\
0 & 0 & 0 & 0 & 1  \\
1 & 0 & 1 & 1 & 0  \\
\end{bmatrix}.$$

The characteristic polynomial of $M_3$ is 
$$x^5-x^3-x^2-1.$$
\end{ex}

The proof of Theorem~\ref{thm:irr} is given in Section~\ref{sec:proof}. Since $x^3-x^2-x-1$ is irreducible mod 3, it is irreducible. The proof for $\Sigma_{5,2}$ and the proof for 
$\Sigma_{2k-1,2}$ are the essentially the same. We give the proof of $\Sigma_{5,2}$ as a warming up.

\begin{ex} \label{ex:2}
$x^5-x^3-x^2-1$ is irreducible.
\end{ex} 

\begin{proof}
$\bullet$ $x^5-x^3-x^2-1$ does not have linear factors.
\begin{enumerate} 
\item[I.] Suppose $x^5-x^3-x^2-1=(x^4+a_3 x^3+a_2 x^2+a_1 x+1)(b_1 x+b_0)$ where $b_0^2=1$ and $b_1^2=1$. 
Then, by comparing the coefficients, we have

\begin{enumerate}
\item[(0)] $b_0=-1$ (the constants), \\
\item[(1)] $b_1+a_1 b_0=0$ (the coefficients of $x^1$), \\
\item[(2)] $a_1 b_1+a_2 b_0=-1$ (the coefficients of $x^2$) $\quad \Rightarrow \quad a_1 b_1+a_2 b_0+1=0$,\\
\item[(3)] $a_2 b_1+a_3 b_0=-1$ (the coefficients of $x^3$) $\quad \Rightarrow \quad a_2 b_1+a_3 b_0+1=0$, \\
\item[(4)] $a_3 b_1+b_0=0$ (the coefficients of $x^4$), \\
\item[(5)] $b_1=1$.
\end{enumerate}

We compute the resultant symmetrically in such a way that when we take the resultant of $(i)$ and $(j)$ with respect to $a_{t}$, we also take the resultant of $(5-i)$ and $(5-j)$ with respect to $a_{4-t}$.

\begin{enumerate}
\item[(6)] $res_{a_1}((1),(2))=-b_1^2+a_2b_0^2+b_0=-1+a_2+b_0$,
\item[(7)] $res_{a_3}((4),(3))=a_2 b_1^2+b_1-b_0^2=a_2+b_1-1$,
\end{enumerate}

Now,

\begin{equation*}
res_{a_2}((6),(7))=b_1-b_0
\end{equation*}

But then, since $b_1=1$, $b_0=1$ which is inconsistent with $(0)$.
\item[II.] Suppose $x^5-x^3-x^2-1=(a_4x^4+a_3 x^3+a_2 x^2+a_1 x+a_0)(x+1)$ where $a_0^2=1$ and $a_4^2=1$. 
Then, by comparing the coefficients, we have

\begin{enumerate}
\item[(0)] $a_0=-1$ (the constants), \\
\item[(1)] $a_1+ a_0=0$ (the coefficients of $x^1$), \\
\item[(2)] $a_1+a_2=-1$ (the coefficients of $x^2$) $\quad \Rightarrow \quad a_1+a_2+1=0$,\\
\item[(3)] $a_2+a_3=-1$ (the coefficients of $x^3$) $\quad \Rightarrow \quad a_2+a_3+1=0$, \\
\item[(4)] $a_3+a_4=0$ (the coefficients of $x^4$), \\
\item[(5)] $a_4=1$.
\end{enumerate}

We compute the resultant symmetrically in such a way that when we take the resultant of $(i)$ and $(j)$ with respect to $a_{t}$, we also take the resultant of $(5-i)$ and $(5-j)$ with respect to $a_{4-t}$.

\begin{enumerate}
\item[(6)] $res_{a_1}((1),(2))=a_2-a_0+1$,
\item[(7)] $res_{a_3}((4),(3))=a_2-a_4+1$.
\end{enumerate}

Now,

\begin{equation*}
res_{a_2}((6),(7))=a_0-a_4
\end{equation*}

But then, since $a_4=1$, $a_0=1$ which is inconsistent with $(0)$.

\end{enumerate} 

$\bullet$ $x^5-x^3-x^2-1$ does not have quadratic factors.
\begin{enumerate} 
\item[I.] Suppose $x^5-x^3-x^2-1=(x^3+a_2 x^2+a_1 x+1)(b_2 x^2+b_1 x+b_0)$ where $b_0^2=1$ and $b_2^2=1$.
Similarly, we symmetrically take the resultant several times until we get two equations with three variables $b_2$,$b_1$, and $b_0$.
Now, taking the resultant one more time with respect to $b_1$, we get
 \begin{equation*}
 4(1-b_0 b_2)=4b_2(b_2-b_0) .
\end{equation*}
But then, the system becomes inconsistent.
\item[II.] Suppose $x^5-x^3-x^2-1=(a_3 x^3+a_2 x^2+a_1 x+a_0)(x^2+b_1 x+1)$.
Similarly, we have $2(1-a_0 a_3)=2a_3(a_3-a_0).$
But then, the system becomes inconsistent.
\end{enumerate}

Therefore, by two bullets, $x^5-x^3-x^2-1$ is irreducible.
\end{proof}

\section{proof of Theorem~\ref{thm:irr}} \label{sec:proof}

We may assume that $k > 3$.

Suppose $$x^{2k-1}-x^k-x^{k-1}-1=  \left(\sum_{i=0}^{l} a_i x^i \right)  \left(\sum_{j=0}^{m} b_j x^j \right)$$ where $l+m=2k-1$.
Without loss of generality, we assume that $l$ is even. We assume that $b_0=b_m=1$  and $a_0^2=a_l^2=1$ (or $a_0=a_l=1$ and $b_0^2=b_l^2=1$). Then, we symmetrically take the resultant several times until we get two equations with three variables, $a_l$, $a_{l/2}$, and $a_0$ ( or $b_m$, $a_{l/2}$, and $b_0$). 
Now taking the resultant one more time with respect to $a_{l/2}$, we have an equation of $a_l$ and $a_0$ (or $b_m$ and $b_0$). We took the resultant so that the coefficients of $a_l$ and $a_0$ (or $b_m$ and $b_0$) are the same η in the final polynomial,  but have different signs. 
If the $\eta \neq 0$, we get an inconsistent system.  At each step, we carefully choose a pair of polynomials so that the two polynomials in the final resultant are not equal up to the powers of $a_l$ and $a_0$ (or $b_m$ and $b_0$). This ensures $\eta \neq 0$. 
Reducing variables through substitution before taking the resultants, greatly reduces work.  Also, the lower the degree, the fewer fakes there are.

\section{proof of the Theorem~\ref{thm:main} on $\Sigma_{2k-1,2}$}

Since $\text{dim}(H_1(N_{2k},\mathbb{R}))=2k-1$, Theorem~\ref{thm:main} is proved by Theorem~\ref{thm:irr}.
\section{acknowledgement}
We thank Erwan Lanneau, Livio Liechti, Julien Marché, Alan Reid, and Darren Long. 
This work was supported by Basic Science Research Program through the National Research Foundation of Korea (NRF) funded by the Ministry of Education, Science and Technology (No. NRF-2018005847). 
\providecommand{\bysame}{\leavevmode\hbox to3em{\hrulefill}\thinspace}
\providecommand{\MR}{\relax\ifhmode\unskip\space\fi MR }
\providecommand{\MRhref}[2]{%
  \href{http://www.ams.org/mathscinet-getitem?mr=#1}{#2}
}
\providecommand{\href}[2]{#2}

\end{document}